\documentclass[conf]{new-aiaa}
 \pdfoutput=1

\usepackage[utf8]{inputenc}

\usepackage{graphicx}
\usepackage{amsmath}
\usepackage[version=4]{mhchem}
\usepackage{siunitx}
\usepackage{longtable,tabularx}
\usepackage{algorithm}
\usepackage{algpseudocode}
\usepackage{tikz}
\usetikzlibrary{shapes.multipart}

\usepackage{xcolor}
\usepackage{glossaries}

\title{Improving Computational Efficiency for Powered Descent Guidance via Transformer-based Tight Constraint Prediction}

\author{Julia Briden\footnote{Doctoral Student, Department of Aeronautics and Astronautics, 77 Massachusetts Avenue, Cambridge, Massachusetts, 02139 USA, and AIAA Student Member.} and Trey Gurga \footnote{Undergraduate Student, Department of Aeronautics and Astronautics, 77 Massachusetts Avenue, Cambridge, Massachusetts, 02139 USA, and AIAA Student Member.}}

\affil{Department of Aeronautics and Astronautics, Massachusetts Institute of Technology, 77 Massachusetts Avenue, Cambridge, Massachusetts, 02139 USA}

\author{Breanna Johnson \footnote{HiMOM Principal Investigator, Flight Mechanics and Trajectory Design Branch, EG5, NASA Johnson Space Center, 2101 E NASA Pky, Houston, TX 77058 USA, and AIAA Atmospheric Flight Mechanics Technical Committee Member.}}

\affil{NASA Johnson Space Center, 2101 East NASA Parkway, Houston, Texas, 77058 USA}

\author{Abhishek Cauligi \footnote{Robotics Technologist, Jet Propulsion Laboratory, California Institute of Technology, Pasadena, CA 91109, USA.}}

\affil{NASA Jet Propulsion Laboratory, California Institute of Technology, 4800 Oak Grove Dr, Pasadena, CA 91109 USA}

\author{Richard Linares \footnote{Rockwell International Career Development Professor and Associate Professor, Department of Aeronautics and Astronautics, 125 Massachusetts Avenue. Senior Member AIAA.}}

\affil{Department of Aeronautics and Astronautics, Massachusetts Institute of Technology, 77 Massachusetts Avenue, Cambridge, Massachusetts, 02139 USA}

\begin{document}

\maketitle

\begin{abstract}
In this work, we present Transformer-based Powered Descent Guidance (T-PDG), a scalable algorithm for reducing the computational complexity of the direct optimization formulation of the spacecraft powered descent guidance problem. T-PDG uses data from prior runs of trajectory optimization algorithms to train a transformer neural network, which accurately predicts the relationship between problem parameters and the globally optimal solution for the powered descent guidance problem. The solution is encoded as the set of tight constraints corresponding to the constrained minimum-cost trajectory and the optimal final time of landing. By leveraging the attention mechanism of transformer neural networks, large sequences of time series data can be accurately predicted when given only the spacecraft state and landing site parameters. When applied to the real problem of Mars powered descent guidance, T-PDG reduces the time for computing the 3 degree of freedom fuel-optimal trajectory, when compared to lossless convexification, from an order of 1-8 seconds to less than 500 milliseconds. A safe and optimal solution is guaranteed by including a feasibility check in T-PDG before returning the final trajectory. 
\end{abstract}

\section{Nomenclature}

\textbf{Variables}{\renewcommand\arraystretch{1.0}
\noindent\begin{longtable*}{@{}l @{\quad=\quad} l@{}}
$\alpha: \mathbb{R}$  & fuel consumption rate equal to $\frac{1}{I_{\text{sp}} g_e}$ [s/m] \\
$\gamma_p: \mathbb{R}$  & maximum allowable tilt angle [rad] \\
$\omega: \mathbb{R}$ & planetary angular velocity [rad/s] \\
$\omega^{\times}: \mathbb{R}^{3 \times 3}$ & skew-symmetric matrix representing the cross product $\omega \times (\cdot)$ \\
$\phi: \mathbb{R}$ & engine angle [rad] \\
$\rho: \mathbb{R}$  & thrust magnitude bound [N] \\
$\xi: \mathbb{R}$ & thrust vector control input slack variable divided by the vehicle's mass [N/kg] \\
$\mathcal{C} \subseteq \mathbb{R}^n$ & feasible set defined by the problem constraints \\
$D(t): \mathbb{R}^{3}$ & drag [N] \\
$d_k$ & number of columns in the transformer weight matrix \\
$g(t): \mathbb{R}^{3}$  & planetary gravitational field [m/s$^2$] \\
$g_e: \mathbb{R}$  & Earth's standard gravitational acceleration [m/s$^2$] \\
$H_{\text{gs}}: \mathbb{R}^{3 \times 3}$  & glideslope constraint matrix \\
$h_{\text{gs}}: \mathbb{R}^{3}$  & glideslope constraint vector \\
$h = 1,...,H$  & transformer head \\
$I_{sp}$ & specific impulse [s] \\
$K_h$ & transformer key matrices \\
$L(t): \mathbb{R}^{3}$ & lift [N] \\
$m: \mathbb{R}$  & spacecraft mass [kg] \\
$N: \mathbb{N}$  & number of discretization nodes \\
$O_h$ & transformer attention output \\
$Q_h$ & transformer key matrices \\
$r(t): \mathbb{R}^{3}$  & spacecraft position [m] \\
$T_c: \mathbb{R}^{3}$  & thrust vector control input [N] \\
$t_0: \mathbb{R}$ & initial time [s] \\
$t_f: \mathbb{R}$ & final time [s] \\
$u: \mathbb{R}^{3}$ & thrust vector control input variable divided by the vehicle's mass [N/kg] \\
$V_h$ & transformer value matrices \\
$v(t): \mathbb{R}^{3}$  & spacecraft velocity [m/s] \\
$z: \mathbb{R}$ & natural log of the vehicle's mass [kg] \\
\end{longtable*}}

\textbf{Functions}{\renewcommand\arraystretch{1.0}
\noindent\begin{longtable*}{@{}l @{\quad=\quad} l@{}}
$\boldsymbol{f}:  \mathbb{R}^{n_x} \times \mathbb{R}^{n_u} \times \mathbb{R} \rightarrow  \mathbb{R}^{n_x}$ & state dynamics function \\
$J: \mathbb{R}^n \rightarrow \mathbb{R}$ & cost function \\
$\boldsymbol{g}: \mathbb{R}^{n_x} \times \mathbb{R}^{n_u} \times \mathbb{R} \rightarrow  \mathbb{R}^{n_c}$ & vector of inequality constraints \\
$\boldsymbol{b}: \mathbb{R}^{n_x} \times \mathbb{R}^{n_x} \times \mathbb{R} \rightarrow  \mathbb{R}^{n_b}$ & vector of boundary constraints \\
$L: \mathbb{R}^{n_x} \times \mathbb{R}^{n_u} \times \mathbb{R} \rightarrow  \mathbb{R}$ & running cost \\
$L_f: \mathbb{R}^{n_x} \times \mathbb{R} \times \mathbb{R} \rightarrow  \mathbb{R}$ & terminal cost \\
$\boldsymbol{u}:[0,t_f] \rightarrow \mathbb{R}^{n_u}$ & control input trajectory \\ 
$\boldsymbol{x}: [0,t_f] \rightarrow \mathbb{R}^{n_x}$ & state trajectory decision variable \\
\end{longtable*}}

\textbf{Notation}{\renewcommand\arraystretch{1.0}
\noindent\begin{longtable*}{@{}l @{\quad=\quad} l@{}}
$(\tau(\theta), t_f^*(\theta))$ & optimal strategy: the set of tight constraints, $\tau(\theta)$, and the optimal final time, $t_f^*(\theta)$ \\
$\hat{e}$  & unit direction vector \\
$\theta$ & parametric problem parameters \\
\end{longtable*}}

\section{Introduction}
\lettrine{N}{ASA's} Moon to Mars Strategy establishes precision landing capabilities as one of the major technologies required for robust exploration and a sustained lunar presence \cite{NASA2023}. To accomplish precision landing in uncertain planetary entry, descent, and landing scenarios, safe and autonomous guidance trajectories must be computed on the order of milliseconds. Currently, precision guidance algorithms fall into two categories: direct and indirect methods. Direct methods discretize the continuous optimization problem and solve it as a parameter optimization problem by numerical optimization (often by primal-dual interior point methods (IPMs)). Indirect methods solve for the necessary conditions of optimality by solving a two-point boundary value problem (TPBVP) corresponding to the state and costate dynamics and their boundary conditions. While indirect methods, including Universal Powered Guidance (UPG) and Propellant-Optimal Powered Descent Guidance, allow for onboard computation of minimum error trajectories, significant simplifications of the spacecraft dynamics and constraints are required to obtain analytical or closed-form solutions; dynamics resulting in a bang-bang profile and linear gravity assumptions are required for problem formulation and inequality constraints on the state and control vectors are not easily implementable in a root-finding framework \cite{Lu2012, Lu2018}. Furthermore, TPBVP-based methods are not guaranteed to find a globally-optimal solution (the problem formulation is not guaranteed to be convex) or a solution in finite time; poor initial guesses may result in convergence failures for gradient-based methods \cite{Lu2018}. On the other hand, direct methods are able to solve generalized powered descent guidance problems at an increased cost due to the cubic time complexity of the solver \cite{Malyuta2021}. Compared to indirect methods, direct methods,  including lossless convexification (LCvx) and successive convexification (SCvx) offer convergence guarantees for generalized problem formulations. LCvx ensures convergence to a globally-optimal solution, while SCvx, under the linear independence constraint qualification (LICQ) and with sufficiently large virtual control penalty weights, provides superlinear convergence to a locally-optimal solution \cite{Acikmese2007ConvexPDG, Malyuta2021}. Additionally, GuSTO, another direct method, leverages information from dual variables to accelerate the algorithm, achieving quadratic convergence rates towards a locally-optimal solution \cite{Bonalli2019GuSTO}. LCvx uses a relaxation of the nonconvex 3 Degree of Freedom (DoF) fuel-optimal powered descent guidance problem with nonconvex control constraints, based on Pontryagin's minimum principle, to obtain a globally optimal solution. SCvx expands upon LCvx by successively convexifying the nonconvex problem to allow for more general 6 DoF problem formulations, including nonconvex state constraints. Unfortunately, this comes at the cost of increased computational complexity. With demonstrated onboard solution times of $2-3$ seconds for LCvx, compared to the radiation-hardened processor requirement of less than 1 second, this work will focus on improving the computational efficiency of direct powered descent guidance algorithms to allow for real-time use of direct optimization methods without the need for custom solver implementations \cite{Acikmese2008, NASA2015}.

Approaches for improving the computational efficiency of onboard direct guidance algorithms include custom solver implementations, storing previously-computed solutions in lookup tables, and informed initial guess strategies; Dueri et al. develop a custom second-order cone-programming (SOCP) IPM in C which exploits the structure of the planetary pinpoint landing powered descent problem to decrease run times by two to three orders of magnitude \cite{Dueri2017}. Since the size of the auto-generated C code grows rapidly with problem size, only problems with less than 1000 solution variables, the length of our decision vector, $\mathbf{x}$, can be solved without the risk of processor instruction cache misses, necessitating slower memory access and reducing computational efficiency \cite{Dueri2017}. Moreover, the use of custom solvers presents a significant challenge for onboard verification due to their size and specificity. As they can span thousands of lines of code, each adjustment or adaptation requires re-verification for flight usage. This can be time-consuming and resource-intensive. In contrast, off-the-shelf solvers, once verified, can more readily adapt to new yet similar problems without the need for constant re-verification. In response to these challenges, a low code footprint first-order solver was developed for the proportional-integral projected gradient (PIPG) method, which can solve the discretized 3 DoF powered descent guidance problem in less than 5 milliseconds for about 25 temporal nodes \cite{Elango2022}. Additionally, a sequential conic optimization (SeCO) solver and a dual quaternion representation, PIPG$_\text{custom}$, were each generated to allow for state-triggered constraints and non-uniform time grids in 6 DoF powered descent guidance problem formulations \cite{kamath2022realtime, kamath2023customized}. Benchmarking results boast 13.7 ms solve times for 16 discretization nodes with SeCO and under 70 ms solve times for 25 discretization nodes with PIPG$_\text{custom}$. While this work significantly improves upon previous custom solvers, this codebase is not immediately generalizable for alternative problem formulations.

An approach which provides more runtime guarantees at the cost of greater memory allocation requirements is the use of interpolable lookup tables. Scharf et al. replace the time-of-flight line search and optimal propellant mass variable with interpolated values from a pre-computed lookup table \cite{Scharf2015}. While this approach enables less than 1 second computation times, interpolation-based lookup tables are not verified for nonconvex problem formulations and these tables can only be applied locally. In addition, time infeasibility can be encountered with too coarse of a sampling grid. An alternative approach to improve the computational efficiency of onboard powered descent guidance includes the development of problem-specific initial guess strategies which predict the optimal state and control trajectory to reduce the number of iterations required to achieve convergence. Kim et al. utilize guided policy search to develop initial guesses for the SCvx algorithm for fixed final time \cite{Kim2022}. The initialization performance for 18 different combinations of penalized trust region (PTR) parameters resulted in a 100\% convergence success rate and a reduction in iteration count when compared to a straight-line initialization. Additionally, Li and Gong apply a feedforward deep neural network (DNN) to predict the optimal final time for LCvx \cite{Li2022}. Use of a DNN resulted in a millisecond-level prediction error but at a cost of 1.2 million training samples. Since the DNN must be retrained on millions of samples for new mission designs and was only shown to be valid for a specified initial state range, this architecture is not conducive to rapid trajectory design.

Rather than naively learning directly the optimal solution for sets of problem parameters, recent work in the field of optimization instead focuses on predicting an interpretable map between parameters and optimal solution strategies \cite{Bertsimas2021, Bertsimas2022}. Since many real-world problems require multiple similar optimization problems to be solved with changing parameters, this method allows data gathered from previous algorithm runs to be used to inform subsequent runs. Bertsimas and Stellato use optimal classification trees (OCTs) and feedforward neural networks (NNs) to solve a multiclass classification problem to predict the strategy, the set of tight constraints (constraints that are equalities at optimality) and integer variable values, for continuous and mixed-integer convex optimization problems. When applied to the motion planning with obstacle-avoidance problem, the machine learning optimizer (MLOPT) implementation resulted in millisecond-level computation times but only had an accuracy close to 10\%, the lowest of all test cases MLOPT was applied to \cite{Bertsimas2022}. While this work proves the feasibility of the NN approach for reducing the computational complexity of online optimization, alternative NN formulations are required to accurately capture time-dependent and nonlinear problem mappings.

When compared to feedforward neural networks (NNs), recurrent neural networks (RNNs)—including a specific type of RNN known as long short-term memory (LSTM) NNs—have been demonstrated to perform well for sequence modeling and time series prediction tasks \cite{Brezak2012, Gers2001}. Still, the sequential computation requirements of recurrent models often suffer from the vanishing or exploding gradient problem, in addition to limitations in parallelization during training. Recently, attention-based mechanisms have been utilized to build transformer NNs which extract global dependencies while allowing for parallelization \cite{Vaswani2017}. Self-attention layers are able to capture these long-range dependencies since their maximum path length that signals the need to traverse is only $O(1)$. Unlike RNNs or LSTMs, transformers process entire sequences in parallel, rather than sequentially. This allows the model to directly learn relationships between elements in the input sequence, regardless of their position. Furthermore, multi-head attention allows the model to pay attention to different parts of the input for feature subspaces, which is imperative for identifying the complicated relationships between problem parameters and the set of tight constraints in a powered descent guidance problem. While widely known for their achievements in natural language processing, including the generative pre-trained transformer (GPT) models and bidirectional encoder representations from transformers (BERT), transformers have also been recently been applied to long-term forecasting applications for time series datasets \cite{Brown2020, Devlin2019, Nie2023, briden2023transformer}. Nie et al. apply a transformer encoder, PatchTST, to multivariate time series prediction and transfer learning for time series benchmarking problems \cite{Nie2023}. Long term forecasting for PatchTST achieves mean squared errors (MSEs) as low as 0.15, exceeding state-of-the-art (SOTA) performance in weather, traffic, and electricity prediction. Additionally, the transformer model maintains low MSE when pre-trained on the Electricity dataset and transferred to other datasets. Therefore, transformer-based NNs efficiently capture highly nonlinear system dynamics in time series data and even generalize patterns to new applications.

While significant progress has been made towards achieving real-time implementable planetary powered descent guidance, several limitations in the SOTA still exist: (i) custom problem-specific solvers or initial guess generators are required to obtain realistic computation times for onboard usage, (ii) data-driven methods offer increased computational efficiency for generalized mixed-integer optimization problems, but they have not yet been shown to work well for long-horizon trajectory generation, and (iii) nonconvex optimization problems, which cannot be losslessly convexified, do not assure a globally-optimal solution and solving these problems often involves multiple iterations of convex approximations, which can be computationally intensive and time-consuming.
    
    \subsection{Contributions}

    \begin{figure}[hbt!]
    \centering
    \includegraphics[width=.9\textwidth]{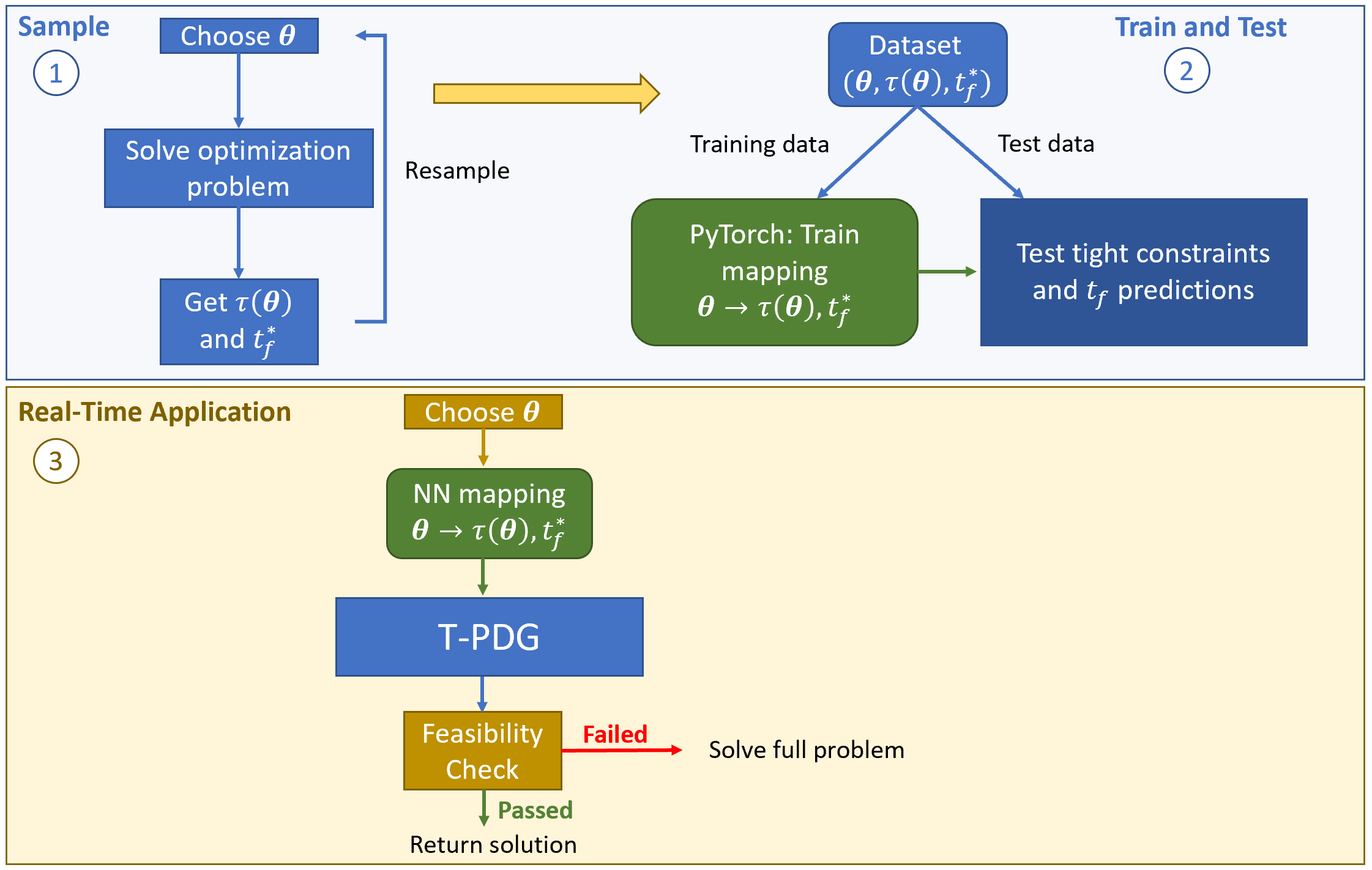}
    \caption{Transformer-based powered descent guidance algorithm overview. Sampling of trajectory datasets and training and testing of the neural networks are shown in steps 1 and 2 (top). The real-time application of T-PDG in a guidance architecture is shown in 3 (bottom).}
    \label{fig:flowchart}
    \end{figure}

    This paper develops what is, to the best of the authors' knowledge, the first transformer-based powered descent guidance algorithm, T-PDG. T-PDG leverages the inherent computational efficiency of transformer neural networks (NNs), which fundamentally operate through efficient matrix multiplication operations. The algorithm is trained to learn a mapping between the input parameters of a powered descent guidance problem and its output: the set of tight constraints and final time at the optimal solution. We denote this output as the problem's ``strategy,'' closely relating to the optimal solution strategies defined in \cite{Bertsimas2021, Bertsimas2022}. Figure \ref{fig:flowchart} illustrates the full architecture for T-PDG.
    
    The design process is segmented into three parts: an offline sampling of problem parameters and their optimal strategies, training and testing of the constraint-prediction NN and final time-prediction NN, and an online real-time application of T-PDG for computationally efficient powered descent guidance. Not only does T-PDG have a low code footprint, when compared to problem-specific solvers, but it also mitigates the typical feedforward NN requirement for a large number of training samples; the entire dataset for training, validation, and test was under 250,000 samples, a significant reduction from the 1.2 million needed for a feedforward NN. T-PDG can be applied to both convex and nonconvex problems, with the only requirement being that some feasible solution exists. When compared to the SOTA, the T-PDG formulation offers a promising balance of computational efficiency, memory efficiency, and generalizability. Given its verifiability across a broad spectrum of problem formulations, T-PDG stands as a potential candidate for implementation on flight-grade processors. The Julia and Python-based software package and test cases included in this paper are available on GitHub as T-PDG\footnote{https://github.com/ARCLab-MIT/T-PDG}.
    
    \subsection{Outline}

    The outline of this paper is as follows: Section \ref{sec: Problem Formulation} defines the general formulation of the powered descent guidance problem along with relevant cost functions, constraints, and assumptions. Section \ref{sec: T-PDG} defines the parameter inputs and strategy outputs, the sampling strategy used for generating training and test data, the transformer model architecture, the process used for training, validation, testing of the model, and finally how the transformer-based architecture leads to an explainable and transferable output for a diverse variety guidance problems. The 3 DoF Mars landing benchmark problem is solved with T-PDG in Section \ref{sec: Application} and compared with LCvx to determine solve times and convergence performance. Finally Sections \ref{sec: Future Directions} and \ref{sec: Conclusion} detail the extension of T-PDG to non-convex problems with feasibility requirements, future work in transfer learning, and concluding remarks on the use of T-PDG for onboard trajectory generation.

\section{Problem Formulation: Powered Descent Guidance}
\label{sec: Problem Formulation}

Trajectory generation problems, in the form of Eq. (\ref{eq: MIP}), and specifically the spacecraft powered descent problem, are formulated as optimal control problems. Where mission objectives determine the cost function and equations of motion, in the form of ordinary differential equations (ODEs), are formulated as constraints, in addition to constraints for state and control requirements. The state $\boldsymbol{x}$ is continuous and $\mathcal{C}$ may be either convex or nonconvex.

\begin{equation}
\begin{aligned}
\min_{\boldsymbol{x}} \quad & J(\boldsymbol{x})\\
\textrm{s.t.} \quad & \boldsymbol{x} \in \mathcal{C}   \\
  \label{eq: MIP}
\end{aligned}
\end{equation}

In the context of the spacecraft guidance problem, Eq. (\ref{eq: MIP}) becomes a semi-infinite optimization problem, as described in Eqs. (\ref{eq: general cost})-(\ref{eq: general boundary constraints}):

\begin{subequations}
\begin{equation}
\min_{t_f,\boldsymbol{u}} \quad L_f (t_0,t_f,\boldsymbol{x}(t_0),\boldsymbol{x}(t_f)) + k \int_{t_0}^{t_f} L(\boldsymbol{x}(\tau),\boldsymbol{u}(\tau),\tau)d \tau 
\label{eq: general cost}
\end{equation}

\begin{equation}
\textrm{s.t.} \quad \dot{\boldsymbol{x}} = \boldsymbol{f}(\boldsymbol{x}(t),\boldsymbol{u}(t),t), \quad \forall t \in [t_0, t_f], 
\label{eq: general dynamics}
\end{equation}

\begin{equation}
\quad \boldsymbol{g}(\boldsymbol{x},\boldsymbol{u},t) \leq 0, \quad \forall t \in [t_0, t_f],
\label{eq: general inequality constraints}
\end{equation}

\begin{equation}
\quad \boldsymbol{b}(\boldsymbol{x}(t_0), \boldsymbol{x}(t_f), t_f) = 0.
\label{eq: general boundary constraints}
\end{equation}
\end{subequations}

The structure of Eqs. (\ref{eq: general cost}) - (\ref{eq: general boundary constraints}) includes ODE constraints (Eq. (\ref{eq: general dynamics})), inequality constraints (Eq. (\ref{eq: general inequality constraints})), and boundary constraints (Eq. (\ref{eq: general boundary constraints})) and is infinite-dimensional; the input trajectory has an infinite number of design parameters due to the continuity of time.  To implement this problem for real-time computation on flight-grade processors, it is formulated as a convex parameter optimization problem by numerical optimization and often solved by primal-dual interior point methods (IPM) or first-order methods.

\subsection{The 3 Degree of Freedom Fuel-Optimal Powered Descent Guidance Problem}\label{problem formulation}

The goal of powered descent guidance (PDG) is to find a sequence of thrust commands that guide the spacecraft from its current state to a desired state. Since the mass of fuel, or wet mass, often represents the majority of the vehicle's mass, and human-level missions come with a significant increase in fuel storage requirements, the objective is often minimizing the fuel usage over time, otherwise known as the fuel-optimal PDG problem. When modeled in 3 DoF, the PDG problem formulation treats the vehicle as a point mass (an assumption which holds when the attitude can be controlled in an inner loop faster than the outer translation control loop). The formulation of the 3 DoF fuel-optimal PDG problem is described in Eqs. (\ref{eq: thrust cost})-(\ref{eq: final boundary constraints}):

\begin{subequations}
\begin{equation}
\min_{t_f,T_c,r,v,m} \quad \int_{0}^{t_f} ||T_c||_2 \mathrm{d} t
\label{eq: thrust cost}
\end{equation}

\begin{equation}
\textrm{s.t.} \quad \dot{r}(t) = v(t), \quad \forall t \in [0, t_f], 
\label{eq: general dynamics r}
\end{equation}

\begin{equation}
\quad \dot{v}(t) = g(t) + \frac{T_c(t) + D(t) + L(t)}{m(t)} - \omega^{\times} \omega^{\times} r(t) - 2 \omega^{\times} v(t), \quad \forall t \in [0, t_f], 
\label{eq: general dynamics v}
\end{equation}

\begin{equation}
\quad \dot{m}(t) = - \alpha ||T_c (t)||_2, \quad \forall t \in [0, t_f], 
\label{eq: general dynamics m}
\end{equation}

\begin{equation}
\quad \rho_{\min} \leq ||T_c(t)||_2 \leq \rho_{\max}, \quad \forall t \in [0, t_f],
\label{eq: thrust bounds}
\end{equation}

\begin{equation}
\quad T_c(t)^T \hat{e}_z \geq ||T_c (t) ||_2 \cos(\gamma_p), \quad \forall t \in [0, t_f],
\label{eq: tilt angle}
\end{equation}

\begin{equation}
\quad H_{\text{gs}} r(t) \leq h_{\text{gs}}, \quad \forall t \in [0, t_f],
\label{eq: affine glideslope}
\end{equation}

\begin{equation}
\quad ||v(t)||_2 \leq v_{\max}, \quad \forall t \in [0, t_f],
\label{eq: max velocity}
\end{equation}

\begin{equation}
\quad m_{\text{dry}} \leq m(t_f),
\label{eq: final mass}
\end{equation}

\begin{equation}
\quad r(0) = r_0, \; v(0) = v_0, \; m(0) = m_{\text{wet}},
\label{eq: initial boundary constraints}
\end{equation}

\begin{equation}
\quad r(t_f) = r_f, \; v(t_f) = 0.
\label{eq: final boundary constraints}
\end{equation}

\end{subequations}

The translational dynamics constraints (Eqs. (\ref{eq: general dynamics r})-(\ref{eq: general dynamics m})) correspond to double integrator dynamics with variable mass viewed in the planet's rotating frame. Eq. (\ref{eq: thrust bounds}) denotes the upper and lower bounds on the rocket engine's thrust vector. The tilt angle constraint, which keeps the spacecraft within $\gamma_p$ of the vertical, is defined in Eq. (\ref{eq: tilt angle}). An affine glideslope constraint is applied in Eq. (\ref{eq: affine glideslope}). Where the glideslope constraint prevents the computed trajectory from going subsurface, the maximum velocity is constrained in Eq. (\ref{eq: max velocity}), and Eq. (\ref{eq: final mass}) constrains the final mass to be greater than or equal to the dry mass such that only the wet mass is used for fuel consumption calculations. Finally, the constraints in Eqs. (\ref{eq: initial boundary constraints})-(\ref{eq: final boundary constraints}) define the initial and final boundary conditions for the spacecraft's state. To formulate this problem as a free final time problem, a line search, or other search method is often used to find a feasible and fuel-optimal final time \cite{Malyuta2021}.

\subsection{The Lossless Convexification (LCvx) Algorithm}\label{sec: LCvx}

 For a full derivation of LCvx for this problem see \cite{Malyuta2021}. The LCvx second-order cone program (SOCP) is defined in Eqs. (\ref{eq: objective function convex})-(\ref{eq: boundary conditions final convex}):

\begin{subequations}
\begin{equation}
\min_{\xi,u,t_f} \int_{0}^{t_f} \xi(t) \mathrm{d} t
\label{eq: objective function convex}
\end{equation}

\begin{equation}
\textrm{s.t.} \quad \dot{r}(t) = v(t), 
\label{eq: dynamics r convex}
\end{equation}

\begin{equation}
\quad \dot{v}(t) = g + u(t) - \omega^{\times} \omega^{\times} r(t) - 2 \omega^{\times} v(t), 
\label{eq: dynamics v convex}
\end{equation}

\begin{equation}
\quad \dot{z}(t) = - \alpha \xi(t), 
\label{eq: dynamics z convex}
\end{equation}

\begin{equation}
\quad \mu_{\min}(t) [1 - \delta z(t) + \frac{1}{2} \delta z(t)^2] \leq \xi(t), 
\label{eq: control constraints lower convex}
\end{equation}

\begin{equation}
\quad \mu_{\max}(t) [1 - \delta z(t)] \geq \xi(t), 
\label{eq: control constraints upper convex}
\end{equation}

\begin{equation}
\quad ||u(t)||_2 \leq \xi(t), 
\label{eq: control constraints norm convex}
\end{equation}

\begin{equation}
\quad u(t)^T \hat{e}_z \geq \xi \cos(\gamma_p), 
\label{eq: control constraints convex}
\end{equation}

\begin{equation}
\quad H_{\text{gs}} r(t) \leq h_{\text{gs}}, 
\label{eq: state constraints glide convex}
\end{equation}

\begin{equation}
\quad ||v(t)||_2 \leq v_{\max}, 
\label{eq: state constraints velocity convex}
\end{equation}

\begin{equation}
\quad \ln(m_{\text{dry}}) \leq z(t_f), 
\label{eq: state constraints dry mass convex}
\end{equation}

\begin{equation}
\quad z_0(t) \leq z(t) \leq \ln(m_{\text{wet}} - \alpha \rho_{\min} t), 
\label{eq: state constraints z bounds convex}
\end{equation}

\begin{equation}
\quad r(0) = r_0, \; v(0) = v_0, \; z(0) = \ln(m_{\text{wet}}), 
\label{eq: boundary conditions initial convex}
\end{equation}

\begin{equation}
\quad r(t_f) = v(t_f) = 0. 
\label{eq: boundary conditions final convex}
\end{equation}
\end{subequations}

Equations (\ref{eq: objective function convex})-(\ref{eq: boundary conditions final convex}) map to the original non-convex formulation through the process of lossless convexification. A slack variable is introduced to remove the nonconvex lower bound in Eq. (\ref{eq: thrust bounds}). Then the variables $\xi$, $u$, and $z$ are used to approximate nonlinear functions of mass. A new objective function, Eq. (\ref{eq: objective function convex}), maximizes final mass, which is equivalent to minimizing fuel consumption. A Taylor series approximation is applied to Eq. (\ref{eq: thrust bounds}) to transform the convex exponential cone constraint into a second order cone. An additional constraint, Eq. (\ref{eq: state constraints z bounds convex}), is then added to ensure the maximum fuel rate is not exceeded.

Since all constraints represent conservative estimates of the original nonconvex constraints, generated solutions will not be infeasible for the original problem.
To ensure the final solution found by the optimization problem is globally optimal the following assumptions must hold:

\begin{enumerate}
    \item For the state $x=\left(r,v\right)\in \mathbb{R}^6$, the state-space matrices,
    $A\ =\ \left[\begin{matrix}0&I_3\\-\omega^\times\omega^\times&-2\omega^\times\\\end{matrix}\right],\ B\ =\ \left[\begin{matrix}0\\I_3\\\end{matrix}\right]$, must be unconditionally controllable.
    
    \item The planet does not rotate about the local vertical of the landing frame ($\omega^\times\widehat{e_z}\neq0$).

    \item The glideslope constraint is only instantaneously active.

    \item For all $\ \theta\in[\frac{\pi}{2}-\gamma_p-\gamma_{\text{gs}},\frac{\pi}{2}+\gamma_p-\gamma_{\text{gs}}]$:

    \begin{enumerate}
        \item $\rho_{\text{min}}\cos\left(\theta\right)<m_{\text{dry}}\left|\left|g\right|\right|_2 \sin\left(\gamma_{\text{gs}}\right)$.

        \item $\rho_{\max}\cos\left(\theta\right)>m_{\text{wet}}\left|\left|g\right|\right|_2 \sin\left(\gamma_{\text{gs}}\right)$.
    \end{enumerate}

    \item The maximum velocity bound is activated at most a discrete number of times.
\end{enumerate}

LCvx is used as the test case for T-PDG in Section \ref{sec: Application} since G-FOLD (Guidance for Fuel Optimal Large Diverts), the algorithm in which LCvx is based upon, is one of the only constraint-satisfying, fuel-optimal, autonomous algorithms that has the potential to scale to to the required 10 km range for powered descent diverts \cite{NASA2015}. 

\section{Transformer-based Powered Descent Guidance (T-PDG)}
\label{sec: T-PDG}

    \subsection{Identifying Optimal Strategies for Solving the Powered Descent Guidance Problem}

    By formulating the powered descent guidance problem into a parametric optimization problem, the parameter vector $\theta$, drawn from a representative set of parameters, $\Theta \subseteq \mathbb{R}^{p}$, is mapped to an optimal set of tight constraints. The tight constraints are defined as the set of constraints that are satisfied as equalities at optimality. If the optimization problem is non-degenerate, the tight constraints serve as support constraints; removal of any tight constraint would result in a decrease in the objective function value for minimization problems \cite{Calafiore}. In both convex and nonconvex optimization problems, the globally optimal solution can occur either on the tight constraint boundaries or within the interior of the feasible set. However, a key difference arises with nonconvex problems, as they may have multiple local optima, each potentially defined by a different set of tight constraints. This complexity requires additional caution: the solutions provided in the parameters and optimal strategy dataset may only be locally optimal and may not necessarily define the global optimum. In this work, we consider only strictly convex problems to eliminate the need for multiple supervision labels for the dataset. While more thought is required to identify a global optimum, or a sufficient local optimum, in nonconvex problems, the methods in this work serve to significantly reduce large optimization problems of either case. Since discretized constrained powered descent guidance problems often have the number of constraints dependent on the number of discretization nodes, methods for tight constraint prediction have the potential to significantly reduce problem size. By defining a reduced problem with only the constraints that the optimal solution pushes against, the optimization algorithm does not have to expend resources checking and managing a large number of constraints that are not critical to finding the optimal solution. In this work, the set of tight constraints is represented as an array of length M, the number of inequality constraints in the constrained optimization problem, where each value in the array is either equal to 1 if the constraint is tight at the optimal solution or 0 otherwise.

    Often time-dependent direct optimization problems have either a fixed-final time or they require an additional embedded optimization problem to determine the optimal, or even a feasible, final time. This is also true for powered descent guidance problems, such as the 3 DoF fuel-optimal guidance problem in Equations (\ref{eq: objective function convex})-(\ref{eq: boundary conditions final convex}) in Section \ref{sec: Problem Formulation} \cite{Acikmese2008, Scharf2015}. As with tight constraints, the optimal final time for an optimization problem can also be formulated as the solution to a parametric optimization problem.

    In the time-dependent optimization problems considered in this study, both the tight constraints and the optimal final time are derived from the same parameter vector, $\theta$, which is drawn from a representative set of parameters, $\Theta \subseteq \mathbb{R}^{p}$. This relationship can be described by the mapping $\theta \rightarrow (\tau(\theta), t_f^*(\theta))$, which we define as the optimal strategy. For the 3 DoF fuel-optimal guidance problem formulated in Equations (\ref{eq: objective function convex})-(\ref{eq: boundary conditions final convex}) in Section \ref{sec: Problem Formulation}, the parameter set is defined as $\theta = \{ \phi, \gamma_{gs}, \gamma_p, r_0, v_0 \}$ (representing engine angle, glideslope angle, pointing angle, initial position, and initial velocity). The strategy for this problem is represented by a binary vector corresponding to the inequality constraints (\ref{eq: control constraints lower convex})-(\ref{eq: state constraints z bounds convex}), $\tau(\theta)$, in addition to a floating point value that represents the final time, $t_f^*$. In the constraint strategy vector, a value of 1 indicates an active constraint, while a value of 0 indicates an inactive constraint at the optimal solution. Given the strategy, a generalized discretized version of Eqs. (\ref{eq: general cost}) - (\ref{eq: general boundary constraints}) are reduced to Eqs. (\ref{eq: general cost new}) - (\ref{eq: general boundary constraints new}):

    \begin{subequations}
    \begin{equation}
    \min_{t_f = t_f^*(\theta), \boldsymbol{u}} \quad L_f (t_0,t_f,\boldsymbol{x}[0],\boldsymbol{x}[N]) + k \sum_{n = 0}^{N} L(\boldsymbol{x}[n],\boldsymbol{u}[n],t[n]) \Delta t 
    \label{eq: general cost new}
    \end{equation}
    
    \begin{equation}
    \textrm{s.t.} \quad \dot{\boldsymbol{x}} = \boldsymbol{f}(\boldsymbol{x}[t],\boldsymbol{u}[t],t), \quad \forall t \in [0, ..., N], 
    \label{eq: general dynamics new}
    \end{equation}
    
    \begin{equation}
    \quad \boldsymbol{g}(\boldsymbol{x},\boldsymbol{u},t) \leq 0, \quad \forall t \in \tau(\theta),
    \label{eq: general inequality constraints new}
    \end{equation}
    
    \begin{equation}
    \quad \boldsymbol{b}(\boldsymbol{x}[0], \boldsymbol{x}[N], t_f) = 0.
    \label{eq: general boundary constraints new}
    \end{equation}
    \end{subequations}

    Where inequality constraints are selectively enforced based on whether they are identified as tight constraints, ensuring computational resources are focused on the constraints which directly influence the optimal solution. Furthermore, the inequality constraints are not changed to equality constraints when identified as tight since this alters the constraint type; for example, a second-order cone (SOC) constraint would generally become nonconvex if changed to an equality constraint. While the original problem may have over 12N+5 constraints (assuming three dimensional equations of motion, one inequality constraint, and five boundary constraints as in Section \ref{sec: LCvx}), the optimal strategy formulation could reduce the number of constraints down to only 3N+5. For the $N = 50$ discretization nodes used in our application problem (Section \ref{sec: Application}), the 605 total constraints could be reduced down to only 155 constraints.
    
    While OCTs and feedforward NNs have been used previously to determine the map between problem parameters and optimal strategies, they were generally less accurate for time-dependent optimal control problems \cite{Bertsimas2022}. Instead we use transformer NNs in this work for strategy prediction, due to their improvements over LSTMs for computationally efficiency and their success in time-dependent forecasting problems \cite{Vaswani2017, Nie2023}. Additionally, for ease of loss function computation and separability, individual mappings for the tight constraints and the optimal final time identification are trained separately by two transformer NNs.

    In the following sections, the sampling strategy used to generate a dataset to train the transformer NNs is discussed, the NN architecture is illustrated, the process used for training, validation, and testing is reviewed, and finally, the interpretability of T-PDG is assessed.

    \subsection{Data Sampling Strategy}

    In order to quickly generate a dataset of relevant parameters and optimal strategy outputs for a general optimization problem, a parallel processing-based sampling algorithm was developed (Algorithm \ref{alg:sampling}). For each variable in $\theta$, a uniform sphere of a selected radius is generated around an initial point. Different radius values are defined for angle parameters, position vectors, and velocity vectors within $\theta$. Points are randomly selected from the sphere to generate a new $\theta$ and this $\theta$ is passed through the optimization solver. The parameter vector $\theta$,the solution's tight constraints, and $t_f$ are stored in an array if the solution is feasible. To determine when there are sufficient samples generated for training, a Good-Turing estimator implementation was explored to generate the probability of unseen strategies yet to be explored in the data set. However, due to the nature of how the tight constraints are generated in a series of zeros and ones, it is very rare to have an exact strategy repeated for trajectories with different initial conditions; just one time step difference in the activation of a tight constraint creates a new strategy. Instead, the algorithm generates a data set of a fixed size and multiple ranges of radii are chosen to sufficiently cover the set of trajectory conditions for training. Based on how large the radii selected are, a set of 100,000 $\theta$ values yields a data set of 50,000 to 80,000 feasible trajectories. Multiple radii are sampled to generate a data set of 300,000 tight constraints for training and testing. The sampling algorithm, programmed in Julia, efficiently leverages parallel computation across 48 cores. As a result, it can process a large dataset of 100,000 data points in just 30 minutes.
    
\begin{algorithm}
  \caption{Uniform Sphere Sampling}\label{alg:sampling}
  \begin{algorithmic}[1]
    \Procedure{Sample}{$\theta_0$, $radius$, $N$=Number of samples}
      \State Initialize an empty list $solutions$
      \For{$i \gets 1$ to $N$} \Comment{Sample $N$ points}
        \State Get a random $\theta$ from a uniform distribution in the range $r$ around $\theta_0$
        \State Generate solution of optimal trajectories of $\theta$
        \State $solution \gets f(\theta)$ \Comment{Function to generate solution tight constraints}

        \If{$solution$ is feasible}
          \State Append $(\theta, solution)$ to $solutions$
        \EndIf
      \EndFor
      \State \textbf{return} $solutions$
    \EndProcedure
  \end{algorithmic}
\end{algorithm}

    \subsection{Transformer Model Structure}

\begin{figure}[hbt!]
    \centering
    \includegraphics[width=0.7\textwidth]{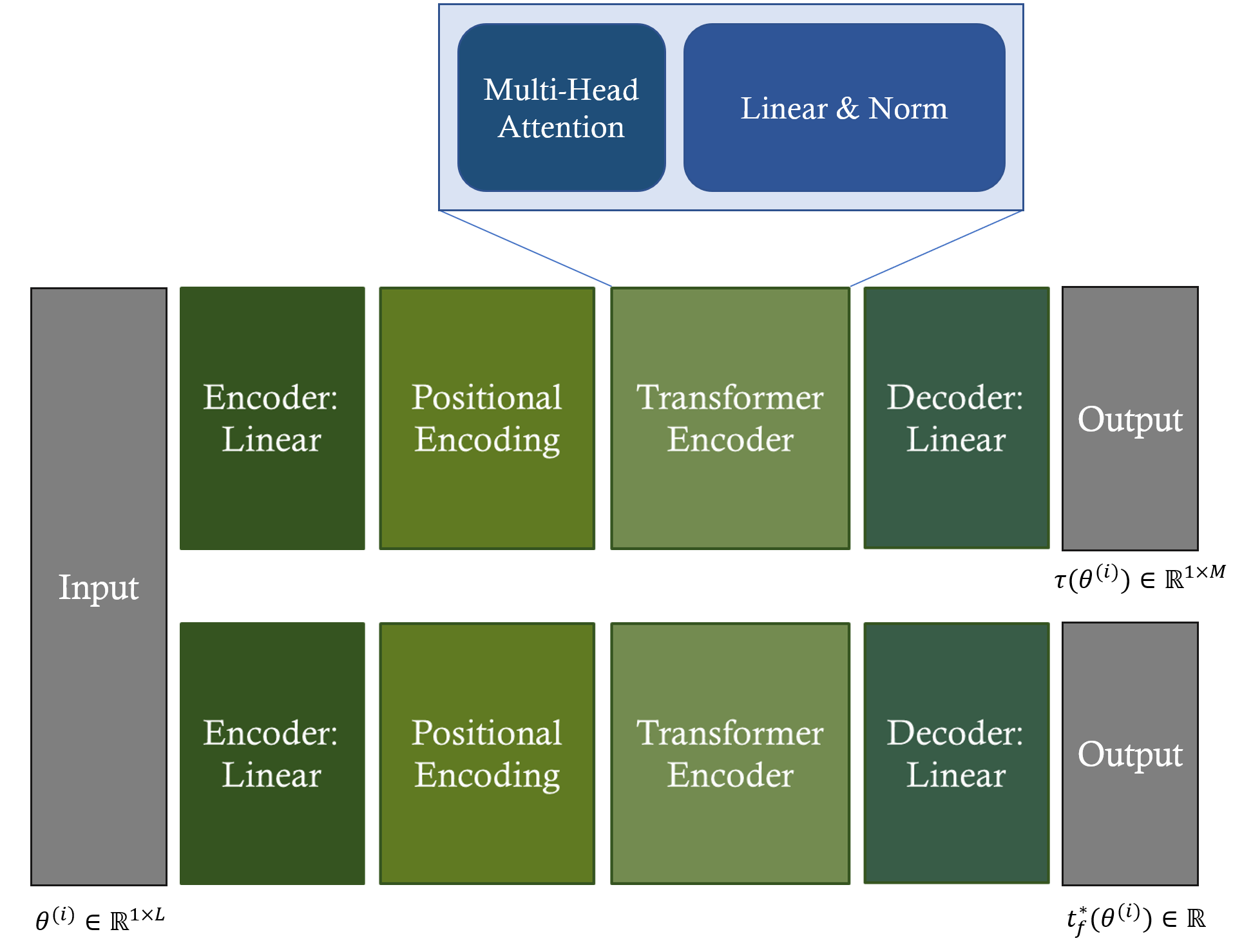}
    \caption{Layers included in the tight constraint prediction and final time transformer neural networks. Both neural networks have the same input, encoder, positional encoder, transformer encoder, and decoder (top and bottom rows). The output of the constraint neural network, $\tau(\theta^{(i)})$, is the set of tight constraints (top right) and the output of the time neural network, $t_f^*(\theta^{(i)})$, is the optimal final time.}
    \label{fig:NN}
    \end{figure}

    The model structure for T-PDG is considered for the following prediction problem: given a set of parametric inputs for a constrained optimization problem, $(\theta_1, ..., \theta_L)$, we would like to predict the set of tight constraints, $\tau(\theta_1), ..., \tau(\theta_L)$, where each $\tau(\theta_i)$ is an 1 x M matrix (M is the number of inequality constraints in the original problem), and the optimal final times, $t_f^*(\theta_1), ..., t_f^*(\theta_L)$. Figure \ref{fig:NN} illustrates the architecture for the constraint prediction and final time prediction NNs.

    Simple linear encoder and decoder layers are used to transfer the input data into a higher dimensional embedding space and the output into a lower dimensional output space. Then a learned position encoding is applied to preserve the temporal order of the input data. From the position encoder, a transformer encoder with number of heads, $h$, uses multi-head attention to transform the data into query matrices, $Q_h^{(i)}$, key matrices, $K_h^{(i)}$, and value matrices, $V_h^{(i)}$. Finally, the attention output is generated by scaled production, as shown in Eq. \ref{eq: Transformer eqn}.

    \begin{equation}
    (O_h^{(i)})^T = \text{Attention} (Q_h^{(i)}, K_h^{(i)}, V_h^{(i)}) = \text{Softmax} (\frac{Q_h^{(i)} K_h^{(i) T}}{\sqrt{d_k}}) V_h^{(i)}
    \label{eq: Transformer eqn}
    \end{equation}

    Additional linear layers, dropout, and LayerNorm layers are also present in the transformer encoder layer. The full model was designed in PyTorch using torch.nn \cite{NEURIPS2019_9015}. The implemented tight constraints and optimal final time NNs for the 3-DoF PDG application (Section \ref{sec: Application}) both have $1 \times 9$-dimensional inputs which include the 3-dimensional initial velocity, 3-dimensional initial position, pointing angle, engine angle, and glideslope angle. Note that the final position and velocity are kept at zero since this application is a powered descent landing problem and reference frames can be adjusted accordingly for a varying final position. Additional parameters for the problem can be included and would only result in a larger input size and a potentially larger required neural network architecture. Furthermore, planetary and spacecraft design parameters were kept constant to represent the chosen mission design. Since state and constraint parameters may change during operation, these variables were chosen as the parameters for the parametric optimization problem.

    \subsection{Real-Time Transformer-based Powered Descent Guidance (T-PDG) Algorithm}

    Algorithm \ref{alg:t-pdg} describes the procedure for applying the transformer NNs for problem reduction in real-time. First, the NN models for predicting tight constraints and optimal final time are called to generate the strategy. Using the strategy, the solver is called to find the corresponding solution and cost, as determined by the problem's cost function. To ensure this returned solution is feasible, it is run as an initial guess on the full-problem solver. Finally, depending on the feasibility of the returned solution, it is either returned or the full problem is solved without the predicted optimal strategy to guarantee that a feasible solution is returned. By implementing this feasibility check, T-PDG is guaranteed to find an optimal solution for any set of parameters which result in a feasible solution for the full problem.

    \begin{algorithm}
      \caption{T-PDG}\label{alg:t-pdg}
      \begin{algorithmic}[1]
        \Procedure{T-PDG}{$\theta$}
            \State $strategy \gets \Call{NN-Prediction}{\mathrm{tight\_constraints\_model}, \mathrm{final\_time\_model}, \theta}$ \Comment{Predict the optimal strategy}
            \State $(soln, cost) \gets \Call{Reduced-Solve}{\theta, strategy}$ \Comment{Solve the strategy-reduced optimization problem}
            \If{$cost \neq \infty$}
                \State $(soln, cost) \gets \Call{Feasibility-Check}{\theta, t_f^*, soln}$ \Comment{Check feasibility of the reduced problem solution}
                \If{$cost = \infty$}
                    \State $(soln, cost) \gets \Call{Full-Solve}{\theta}$ \Comment{Solve the full optimization problem if infeasible}
                \EndIf
            \Else
                \State $(soln, cost) \gets \Call{Full-Solve}{\theta}$
            \EndIf
            \State \textbf{return} $soln$
        \EndProcedure
      \end{algorithmic}
    \end{algorithm}

    \subsection{Training, Validation, and Testing}

    After completing the sampling step, all sampled data is split into 80\% training and validation data and 20\% test data for both constraints prediction and optimal final time prediction. Then mean and standard deviation of the training sets are computed and the train and test datasets are standardized by subtracting the mean and dividing by the standard deviation. For the training process, K-fold cross-validation, with K=2, was employed to split the training processing into training and validation steps. In each iteration of the validation process, one part is used for training the model and the other for testing its performance. This approach enhances the network's ability to generalize to new, unseen data, compared to using a fixed single train-test split. Mean squared error (MSE) loss was used for the training loss function for both NNs and a learning rate scheduler was employed for the final time prediction NN to improve training convergence.

    \subsection{Interpretability and Generalizability}

    As noted in Bertsimas and Stellato, the identification and prediction of optimal strategies for constrained optimization problems not only allow for efficient solution recovery but also provide interpretability and intuition behind the optimal solutions; researchers are able to better understand how problem parameters affect optimal control decisions \cite{Bertsimas2021}. Specifically for transformer NNs, the attention mechanism present in the transformer encoder provides a source of interpretability for the embedding space.
    
    As an analogy, consider the domain of explicit Model Predictive Control (MPC). In explicit MPC, the control problem is solved offline, and the solution space is divided into different regions or 'polytopes'. Each polytope corresponds to a specific set of conditions or parameters within the system. Within each of these polytopes, a specific control law is activated based on the current state of the system. This separation into polytopes serves to provide a clear, interpretable map of how different parameter states lead to different control actions. Comparatively, for transformer NNs, the attention mechanism present in the transformer encoder provides a source of interpretability for the embedding space. By employing dimensionality reduction techniques like t-Distributed Stochastic Neighbor Embedding (t-SNE), high-dimensional representations learned by transformers can be projected down into two or three-dimensional spaces \cite{Maaten2008}. This process is analogous to visualizing the MPC's polytope activations, where each cluster or pattern in the t-SNE plot can be viewed as a 'region' or 'polytope' of operation with similar characteristics in terms of the number of active constraints or the optimal final time value. Furthermore, the generalizability of the T-PDG algorithm can be assessed by comparing the t-SNE embeddings of both the training and testing datasets. If the structure of these embeddings is similar, it indicates effective learning and application of the algorithm, consistent to how a well-generalized MPC controller would frequently apply the correct control laws across various regions in its operational space.
    
    In the field of powered descent guidance, Malyuta et al. suggest that a deeper undiscovered underlying theory for how and why problems can be losslessly convexified exists due to the fact that arbitrarily small perturbations of the dynamics can recover LCvx \cite{Malyuta2021}. The development of methods which result in transparent parameter to optimal strategy mappings enable the potential discovery of this theory, as well as the design of a wider range of generalizable global optimization algorithms. 

\section{Application: 3 Degrees of Freedom T-PDG}
\label{sec: Application}

    \subsection{Problem Setup and Parameters}

    To analyze the performance of T-PDG, the algorithm was trained on data from the LCvx Mars rocket landing problem implementation for 50 discretization nodes, $N = 50$, from the SCP Toolbox, which uses the ECOS solver\footnote{https://github.com/jump-dev/ECOS.jl} in Julia \cite{Malyuta2021}. Since well-known results for this test case are already recorded, they can be easily used to determine the effect T-PDG has on both convergence and efficiency. Table \ref{tab:algorithm_settings} shows the planetary and spacecraft parameters as described in Malyuta et al. \cite{Malyuta2021}.

    The sampled trajectory and the range of values included in the training, validation, and test datasets are detailed in Table \ref{tab:sampled}. The entire dataset comprises 242,293 samples, which encompass the sampled problem parameters, tight constraints, and optimal final times. These parameters are designer-specified and, in this application, they are chosen to be $\theta = { \phi, \gamma_{gs}, \gamma_p, r_0, v_0 }$, representing the engine angle, glideslope angle, pointing angle, initial position, and initial velocity, respectively. These parameters are particularly relevant as they may vary between the parachute phase and the powered descent initiation. The final state is maintained constant, as a change in the reference frame can be used to represent any desired configuration. Furthermore, the strategy for the Transformer-based Powered Descent Guidance (T-PDG) is depicted as a binary vector $\tau(\theta)$, corresponding to the inequality constraints (\ref{eq: control constraints lower convex})-(\ref{eq: state constraints z bounds convex}). This vector is complemented by a floating-point value representing the optimal final time, $t_f^*$.

    \begin{table}
        \caption{\label{tab:algorithm_settings} Algorithm Settings}
        \centering
        \begin{tabular}{ll}
            \hline
            \textbf{Planetary Parameters} & \\
            \hline
            \textbf{Parameter} & \textbf{Value} \\
            Final position ($r_f$) &$0 \hat{e}_x + 0 \hat{e}_y + 0 \hat{e}_z$ m \\
            Final velocity ($v_f$) &$0 \hat{e}_x + 0 \hat{e}_y + 0 \hat{e}_z$ m/s \\
            Gravity ($g$) & $-3.7114 \hat{e}_z$ m/s$^2$ \\
            Latitude of landing site ($\theta$) & $30$ deg \\
            Mars sidereal ($T_{\text{sidereal\_mars}}$) & 88,642.44 s \\
            Standard gravity ($g_e$) & 9.807 m/s$^2$ \\
            \hline
            
            \textbf{Spacecraft Parameters} & \\
            \hline
            \textbf{Parameter} & \textbf{Value} \\
            Dry Mass ($m_{\text{dry}}$) & 1,505 kg \\
            Wet Mass ($m_{\text{wet}}$) & 1,905 kg \\
            Isp & 225.0 s \\
            Number of engines ($n_{\text{eng}}$) & 6 \\
            Max physical thrust of single engine ($T_{\max}$) & 3,100 N \\
            Max allowed velocity ($v_{\max}$) & 138.89 m/s \\
            Final velocity ($v_f$) & $0 \hat{e}_x + 0 \hat{e}_y + 0 \hat{e}_z$ m/s \\
            Min allowed thrust of single engine ($T_1$) & $930$ N \\
            Max allowed thrust of single engine ($T_2$) & $2,480$ N \\
            Min allowed thrust of all engines ($\rho_{\text{min}}$) & $4,971.8$ N \\
            Max allowed thrust of all engines ($\rho_{\text{max}}$) & $13,258.2$ N \\
            Fuel consumption rate ($\alpha$) & $0.0005$ s/m \\
            Number of discretization points ($N$) & 50 \\
            \hline
        \end{tabular}
    \end{table}

    \begin{table}
        \caption{\label{tab:sampled} Sampled Dataset}
        \centering
        \begin{tabular}{lll}
            \hline
            \textbf{Parameter} & \textbf{Sampled Trajectory} & \textbf{Range} \\
            Engine angle ($\phi$) & 10 deg & [0 , 20] deg \\
            Glideslope angle ($\gamma_{gs}$) & 80 deg & [70 , 90] deg \\
            Pointing angle ($\gamma_p$) & 60 deg & [50 , 70] deg \\
            Initial Position ($r_0$) & $2000 \hat{e}_x + 2000 \hat{e}_y + 1000 \hat{e}_z$ m & [1500 , 2500]$\hat{e}_x$ + [1500 , 2500]$\hat{e}_y$ + [500 , 1500]$\hat{e}_z$ m\\
            Initial Velocity ($v_0$) & $-15 \hat{e}_x + -15 \hat{e}_y + -30 \hat{e}_z$ m/s & [-115 , 85]$\hat{e}_x$ + [-115 , 85]$\hat{e}_y$ + [-101.7 , 70]$\hat{e}_z$ m/s\\
            \hline
        \end{tabular}
    \end{table}

    The NNs for tight constraint prediction and final time prediction both used the transformer architecture shown in Figure \ref{fig:NN} from Section \ref{sec: T-PDG}. Since predicting a $1 \times 8N-3$-dimensional array for constraint prediction is a larger problem than predicting a scalar final time, the constraint prediction NN utilized 384-dimensional layers, two heads, and four layers. The final time NN was designed to be much smaller, with 64-dimensional layers, one head, and two layers. Both NNs had a dropout rate of 0.1 to reduce overfitting to the training data.

    \subsection{Results and Analysis}

    For training the constraint prediction NN, the learning rate was set to a constant 0.001. Since convergence for training of the time prediction NN was more parameter-sensitive, the initial learning rate was set at 0.01 and was divided by 4000 until 4000 warm-up steps were reached. After which, $0.01$( current\_step - warmup\_steps$^{-0.5})$ is used for the learning rate. Only one epoch was used for training since both the models converge in less than one epoch. Less than 130,000 samples, at a batch size of 320, were required for the optimal final time prediction NN training to converge. Additionally, only 10,000 samples, at a batch size of 128, were required to train the tight constraint prediction NN. Table \ref{tab:trainandtest} includes the training and validation MSE at the second fold of training, as well as the accuracy of each NN on the test dataset. Due to the the different formats of the tight constraints output and the final time output, MSE is used for final time prediction accuracy and binary accuracy, determined by the number of correct labels divided by the total number of labels, is used for the tight constraints NN.

\begin{table}[ht]
    \caption{Training and Testing of T-PDG}
    \label{tab:trainandtest}
    \centering
    \begin{tabular}{lcccc}
    \hline
    Model & Train (MSE) & Validation (MSE) & Test (MSE / Binary Accuracy) & Number of Parameters \\
    \hline
    Constraint NN & 0.063 & 0.049 & 94.01\% (Binary Accuracy) & 8,829,453 \\
    Final Time NN & 121.14 & 106.55 & 104.26 (MSE) & 563,009 \\
    \hline
    \end{tabular}
\end{table}

    The results obtained by applying T-PDG, from Algorithm \ref{alg:t-pdg}, to LCvx, using 775 samples from the test dataset, are shown in Table \ref{tab:performance}. Both the computation time and feasibility metrics were analyzed for each portion of the T-PDG algorithm and the LCvx algorithm. Note that even for a less than 100\% feasibility in the feasibility check, T-PDG is still 100\% feasible since the full-problem is solved when this occurs.

\begin{table}
    \caption{\label{tab:performance} Performance Comparison of T-PDG with LCvx}
    \centering
    \begin{tabular}{lccc}
    \hline
    Algorithm & Feasibility & \multicolumn{2}{c}{Computation Time [ms]} \\
    & & Mean & Standard Deviation \\
    \hline
    T-PDG NN prediction & \textbf{100 \%} & 6.261 & 2.910 \\
    T-PDG reduced-solve & \textbf{100 \%} & 33.27 & 30.81 \\
    T-PDG feasibility check & 83.74 \% & 333.72 & 634.76 \\
    \hline
    T-PDG Total & \textbf{100 \%} & \textbf{373.25} & 635.40 \\
    \hline
    LCvx & \textbf{100 \%} & 1,673.3 & 352.23 \\
    \hline
    \end{tabular}
\end{table}

\begin{figure}[hbt!]
        \centering
        \begin{minipage}{.5\textwidth}
            \centering
            \includegraphics[width=\linewidth]{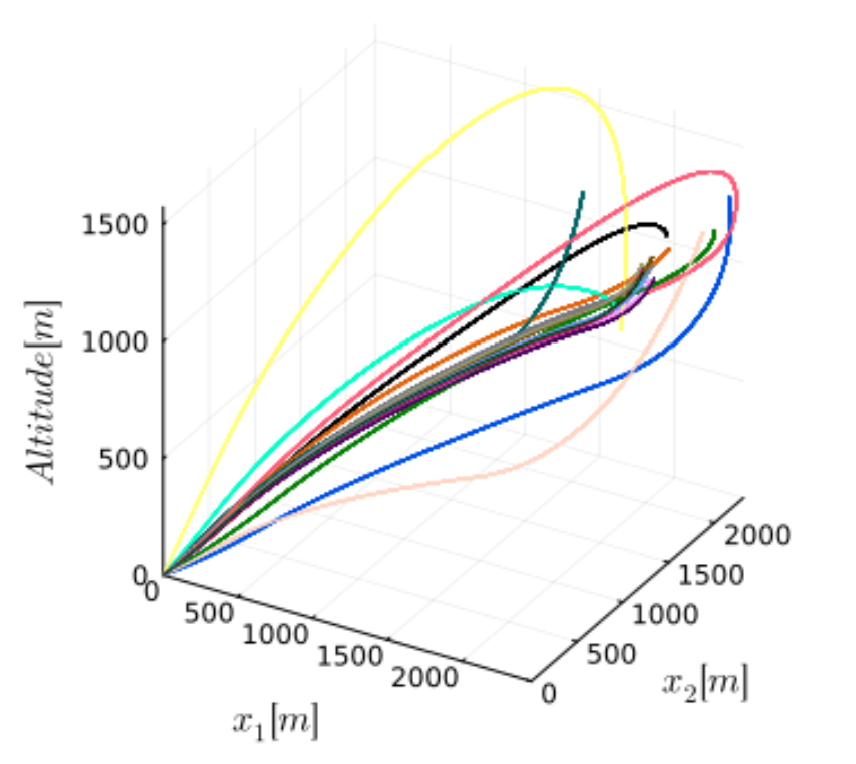}
        \end{minipage}\hfill
        \begin{minipage}{.5\textwidth}
            \centering
            \includegraphics[width=\linewidth]{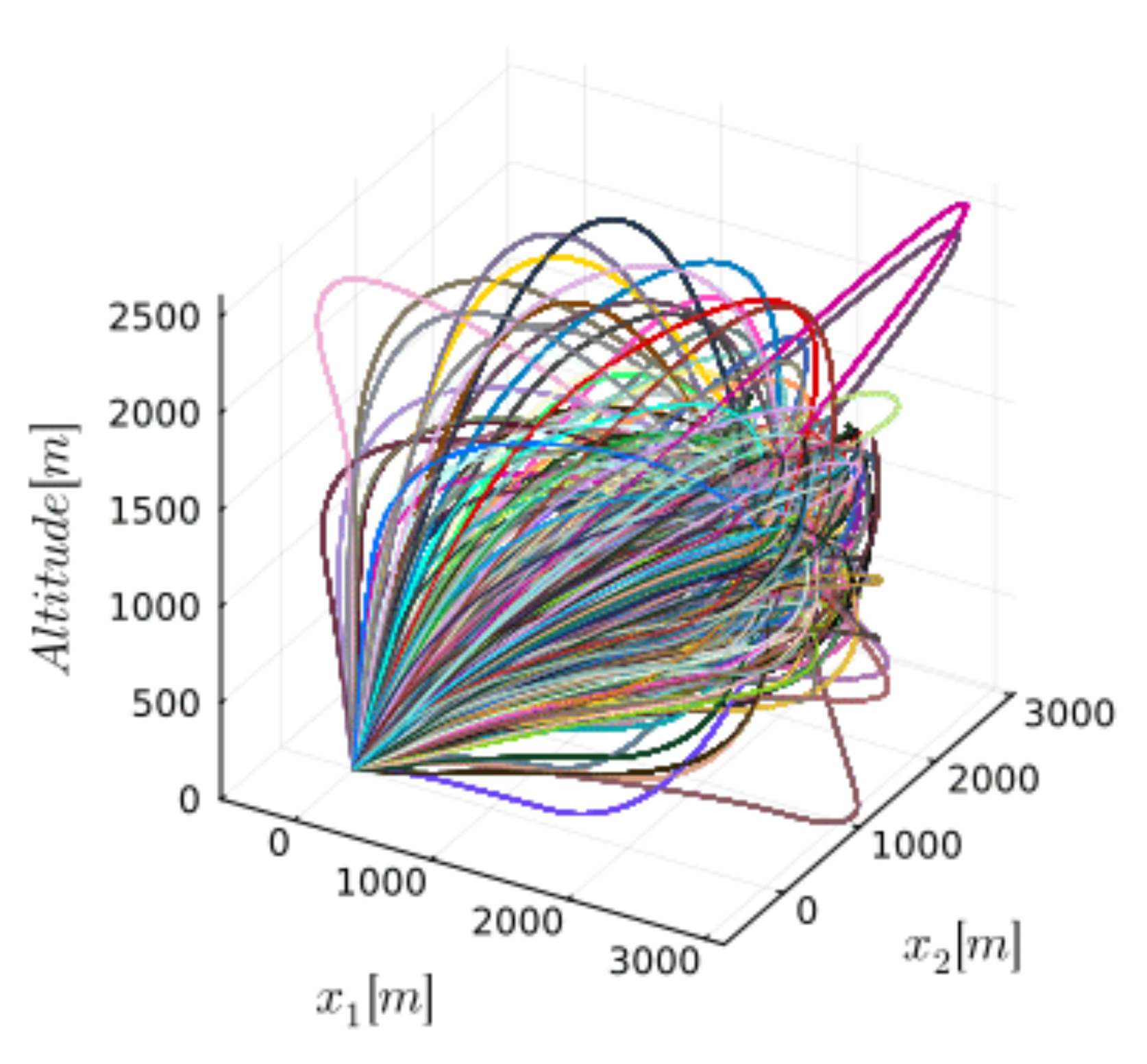}
        \end{minipage}
    \caption{\href{https://youtu.be/1HhaewVzo5Q}{Test trajectories for T-PDG}. Left: first few trajectories obtained by T-PDG. Right: trajectories for the full 775 samples from the test dataset, computed by T-PDG. All constraints are consistently satisfied for each trajectory, evident by exclusively above surface trajectories, maintaining a glideslope angle within the 70-90 degrees range. Furthermore, these trajectories demonstrate dynamics that align with the fuel-minimization objective, while also adhering to the specified engine angle and pointing angle constraints.}
        \label{fig:trajectories}
    \end{figure}

    Figure \ref{fig:trajectories} shows the sets of T-PDG-computed test trajectories as they are first generated (the lefthand side), as well as the full set of 775 test trajectories (the righthand side).

    When compared to LCvx, T-PDG reduces the mean computation time required to compute the fuel-optimal powered descent solution by 78\%. With a mean computation time exceeding 1.5 seconds for LCvx, T-PDG consistently achieves computation times of less than 1 second within one standard deviation. Therefore, preliminary analysis shows T-PDG consistently meeting the <1 second computation time requirement of flight grade processors \cite{NASA2015}. Since this time includes the 16\% of test cases when the feasibility check is not passed and the full solution must be computed, T-PDG presents both an efficient and reliable approach for improving the computational efficiency of the powered descent guidance problem. Instances of infeasibility in T-PDG typically arise from underestimating the optimal final time for the constraint-satisfying trajectory. One possible improvement to enhance initial feasibility could involve incorporating a conservative margin in the final time prediction, though this might slightly increase fuel consumption. Additionally, when compared to the 1.2 million training samples required for the feedforward implementation of final time prediction, less than 10,000 samples were required for the final time prediction NN to converge in T-PDG \cite{Li2022}. From this observation and the convergence of the larger constraint prediction NN after less than 130,000 samples, T-PDG appears to be quite sample efficient. Importantly, the tight constraints prediction problem, as applied to optimal control, is well-predicted by transformer NNs; the constraint prediction NN achieved over 94\% accuracy when evaluated on test data for 50 time discretization nodes. T-PDG and transformer-based data-driven optimization methods provide substantial promise in the prediction of highly interpretable solutions for time-dependent optimization problems \cite{Bertsimas2022}.

    \subsection{Assessment of Interpretability and Generalizability}

\begin{figure}[hbt!]
        \centering
        \begin{minipage}{.5\textwidth}
            \centering
            \includegraphics[width=\linewidth]{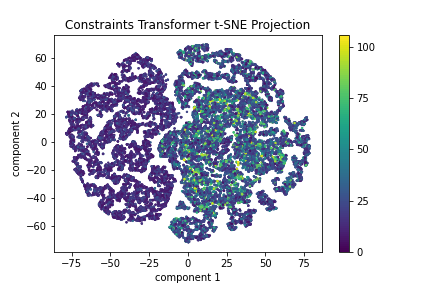}
            \captionsetup{font=small}
            \caption{t-SNE for tight constraints prediction.}
            \label{fig:tsne_constraints_train}
        \end{minipage}\hfill
        \begin{minipage}{.5\textwidth}
            \centering
            \includegraphics[width=\linewidth]{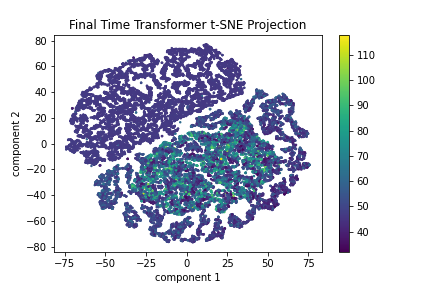}
            \captionsetup{font=small}
            \caption{t-SNE for optimal final time prediction.}
            \label{fig:tsne_time_train}
        \end{minipage}
        \caption{t-SNE visualizations for training dataset. Left: Tight constraints prediction, where each color denotes the number of active constraints. Right: Optimal final time prediction, where each color denotes the value of the optimal final time in seconds. The results of the t-SNE visualization allow for a practitioner to discern relationships between the optimal control problem parameters and tight constraints. In this case, we see that there are two distinct clusters (low vs. high number of active constraints and optimal final times) that emerge for the PDG problem.}
        \label{fig:tsne_training}
    \end{figure}

\begin{figure}[hbt!]
        \centering
        \begin{minipage}{.5\textwidth}
            \centering
            \includegraphics[width=\linewidth]{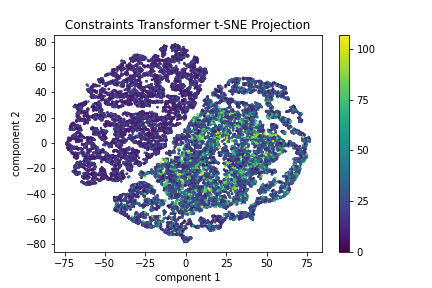}
            \captionsetup{font=small}
            \caption{t-SNE for tight constraints prediction.}
            \label{fig:tsne_constraints}
        \end{minipage}\hfill
        \begin{minipage}{.5\textwidth}
            \centering
            \includegraphics[width=\linewidth]{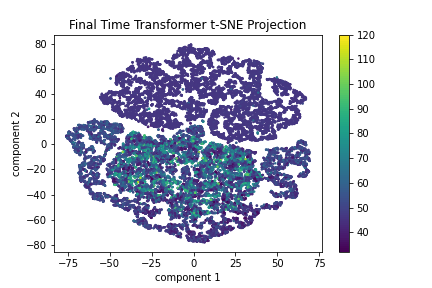}
            \captionsetup{font=small}
            \caption{t-SNE for optimal final time prediction.}
            \label{fig:tsne_time}
        \end{minipage}
        \caption{t-SNE visualizations for test dataset. Left: Tight constraints prediction, where each color denotes the number of active constraints. Right: Optimal final time prediction, where each color denotes the value of the optimal final time in seconds. The results of the t-SNE visualization allow for a practitioner to discern generalizability of the NN to test data, outside of the training dataset. In this case, we see that the two distinct clusters (low vs. high number of active constraints and optimal final times) that emerged in the training data are also apparent for the test dataset, indicating a generalized classification of problem types.}
        \label{fig:tsne_test}
    \end{figure}

    To assess the interpretability and generalizability of each transformer model, 2D t-SNE visualization were generated using sklearn.manifold.TSNE and analyzed for 30,000 training and 30,000 test samples, each with batch size 1 \cite{scikit-learn}. The t-SNE technique, similar to visualizing explicit MPC polytope activations, reduces high-dimensional data to a more interpretable two-dimensional space. This reduction allows for an analysis of how input data cluster in relation to transformer model embeddings. The t-SNE visualization for the training dataset are shown in Figure \ref{fig:tsne_training} and the t-SNE visualization for the test dataset are shown in Figure \ref{fig:tsne_test}. In these figures, 'component 1' and 'component 2' denote the two-dimensional reduced components resulting from the t-SNE process. A comparison between the training and test visualizations provides an understanding of how well the model generalizes to new, unseen data, indicated by similar clustering patterns in both sets of embeddings.

    The color scale in the tight constraints prediction NN represents the number of tight constraints at the optimal solution. This choice of color scale, instead of individual constraint representations, is due to the high discretization of the problem, which results in a limited repetition of strategies. The number of tight constraints still provides a general view of the optimal strategy based on how constrained each parametric problem is. The active constraints occurring most often, based on the sampled dataset, are the second-order cone constraint on thrust (Eqn. \ref{eq: control constraints norm convex}) and the maximum and minimum thrust bounds (Eqns. \ref{eq: control constraints lower convex}-\ref{eq: control constraints upper convex}). The mass constraints (Eqns. \ref{eq: state constraints dry mass convex}-\ref{eq: state constraints z bounds convex}) and pointing angle constraint (Eqn. \ref{eq: control constraints convex}) are the next most often occurring tight constraint, often occurring at the start or end of a trajectory. For the optimal final time NN, the color scale corresponds to the rounded optimal final time value.

   In observing both the training (Figure \ref{fig:tsne_training}) and test datasets (Figure \ref{fig:tsne_test}), similar structures are evident not only between these datasets but also between the tight constraints and final time transformers. Each t-SNE projection displays a two-lobe formation, indicating distinct or predominant divisions within the dataset. In the tight constraints NN, one lobe is characterized by a near absence of active constraints, depicted in dark purple. In the final time NN, the corresponding dark purple lobe in the final time NN represents data with very short optimal final times (less than 55 seconds). These lobes are distinctly separated from other data regions in the embedding space, highlighting a clear classification by both NNs of problems with fewer constraints and shorter final times. The second lobe in the tight constraints NN predominantly features medium to highly constrained samples, shown in shades of blue, green, and yellow. These samples appear to be regularly spaced in both the training dataset and test dataset. These samples are regularly spaced in both the training and test datasets, suggesting that medium to highly-constrained problems are not distinctly separated in the embedding space. However, this pattern's consistency across both the training and test datasets implies that this lack of distinct separation may not adversely impact the generalizability of the tight constraints NN. Since the coloring only indicates the total number of active constraints and not the type or exact values of the active constraints, this other embedding space lobe cannot be fully interpreted in terms of its separability for medium to highly constrained problems. In the final time NN, the second lobe exhibits more defined structures. Here, samples with longer optimal final times (light green to yellow) cluster towards the center, while those with medium to lower final times (dark purple to blue) are found around the lobe's edges. Since this phenomenon occurs in both the training and test datasets, generalizability of the final time NN remains likely high.
   
   Overall, the t-SNE visualizations demonstrate matching structures in the embedding spaces of training and test datasets. Notably, the dual-lobe structure is a common feature in both the tight constraints and optimal final time NNs, suggesting a coherent and interpretable pattern in the model's learning and prediction process.

\section{Future Directions}
\label{sec: Future Directions}

While the results obtained for T-PDG are promising for improving both computational efficiency and interpretability of the powered descent guidance problem, several future directions can be explored further:

\begin{itemize}
\item \textbf{Transfer learning:} We believe that there is significant potential in using a pre-trained transformer on a simpler convex problem, such as LCvx, then fine-tuning it for more complex nonconvex optimization problems. This transfer learning approach could boost the learning efficiency and performance of T-PDG, especially when the training data for the complex problem is scarce or hard to collect. Furthermore, a fine-tuned T-PDG may be used for new mission designs or changes to the chosen set of parameters for the parametric optimization problem. Previous applications of pre-trained transformer architectures in long-term forecasting have been shown to maintain low MSE when transferred between other datasets \cite{Nie2023}. This indicates T-PDG's potential to adapt well to transfer learning approaches to further improve generalizability.
\item \textbf{Multi-task learning:} While the constraints prediction and final time prediction NNs were trained separately to ensure that the loss function was not dominated by one task, they could potentially be trained together in a multi-task learning setup. This could enforce a shared representation learning which could potentially increase the performance of each task. From the t-SNE plots in Figure \ref{fig:tsne_training} and Figure \ref{fig:tsne_test}, similar structures appear in the embedding space for both NNs, indicating a similar representation for both tasks.
\item \textbf{Active learning:} An active learning approach could be adopted to iteratively add more complex samples to the training set. This can help in learning a more robust and generalized model. Since multiple similar powered descent guidance problems are often solved during analysis and operation, active learning frameworks provide the potential to incorporate these new samples into the current learned representation.
\item \textbf{Extended application:} The current study only considers a LCvx Mars rocket landing problem, but T-PDG could be applied to other powered descent guidance problems, as well as to ascents, aborts, and nonconvex optimization problems.
\item \textbf{In-depth interpretability study:} Although we offer an interpretability analysis based on t-SNE visualization, further t-SNE analyses could be conducted to better understand the type of constraints activated within each cluster. Discoveries in explainable tight constraint prediction may also contribute to uncovering foundational theories in optimization for powered descent guidance.
\end{itemize}

\section{Conclusion}
\label{sec: Conclusion}

This work presents T-PDG, a transformer-based approach to predict the solution of LCvx problems for real-time powered descent guidance. By leveraging the attention mechanism of transformer neural networks, T-PDG enables the millisecond-level prediction of optimal strategies which reduce the average runtime of the LCvx from over 1.5 seconds to under 500 milliseconds. Moreover, T-PDG ensures the feasibility of the final solution through a feasibility check and worst-case full solve. T-PDG was applied to a Mars rocket landing problem, offering a 78\% reduction in computation time, when compared to the LCvx algorithm, while retaining 100\% feasibility. In addition, the transformer-based approach exhibits promising interpretability and generalizability, as illustrated by the t-SNE visualizations of the models' embedding space. The algorithm is generalizable to any optimization problem in which a feasible solution exists and no custom codes were required to improve computational efficiency. T-PDG represents a near-term implementable algorithm capable of enabling real-time trajectory generation on flight-grade processors. Future work will focus on extending the application of T-PDG to more complex scenarios, including nonconvex problems in 6 DoF with nonlinear drag and gravity terms.

\section*{Acknowledgments}
This work was supported in part by a NASA Space Technology Graduate Research Opportunity 80NSSC21K1301 and the National Science Foundation under award NSF-PHY-2028125.
The research was carried out at the Jet Propulsion Laboratory, California Institute of Technology, under a contract with the National Aeronautics and Space Administration (80NM0018D0004).
The authors would like to thank the MIT SuperCloud for providing HPC, database, and consultation resources that have contributed to the research results reported in this paper.

\bibliography{sample}

\begin{thebibliography}{28}
\newcommand{\enquote}[1]{``#1''}
\providecommand{\natexlab}[1]{#1}
\providecommand{\url}[1]{\texttt{#1}}
\providecommand{\urlprefix}{URL }
\expandafter\ifx\csname urlstyle\endcsname\relax
  \providecommand{\doi}[1]{\discretionary{}{}{}https://doi.org/#1}\else
  \providecommand{\doi}[1]{\discretionary{}{}{}\urlstyle{rm}\url{https://doi.org/#1}}\fi

\bibitem[{NAS(2023)}]{NASA2023}
\enquote{NASA’s Moon to Mars Strategy and Objectives Development,} 2023.

\bibitem[{Lu(2012)}]{Lu2012}
Lu, P., \enquote{A Versatile Powered Guidance Algorithm,} \emph{AIAA Guidance,
  Navigation, and Control Conference}, 2012.
\newblock \doi{10.2514/6.2012-4843}.

\bibitem[{Lu(2018)}]{Lu2018}
Lu, P., \enquote{Propellant-Optimal Powered Descent Guidance,} \emph{JGCD},
  Vol.~41, No.~4, 2018, pp. 813, 826.
\newblock \doi{10.2514/1.G003243}.

\bibitem[{Malyuta et~al.(2021)Malyuta, Reynolds, Szmuk, Lew, Bonalli, Pavone,
  and Acikmese}]{Malyuta2021}
Malyuta, D., Reynolds, T., Szmuk, M., Lew, T., Bonalli, R., Pavone, M., and
  Acikmese, B., \enquote{Convex Optimization for Trajectory Generation,}
  \emph{IEEE Control Systems Magazine}, Vol.~42, No.~5, 2021, pp. 40--113.

\bibitem[{Açıkmeşe and Ploen(2007)}]{Acikmese2007ConvexPDG}
Açıkmeşe, B., and Ploen, S.~R., \enquote{Convex Programming Approach to
  Powered Descent Guidance for Mars Landing,} \emph{Journal of Guidance,
  Control, and Dynamics}, Vol.~30, No.~5, 2007.
\newblock \doi{10.2514/1.27553}.

\bibitem[{Bonalli et~al.(2019)Bonalli, Cauligi, Bylard, and
  Pavone}]{Bonalli2019GuSTO}
Bonalli, R., Cauligi, A., Bylard, A., and Pavone, M., \enquote{{GuSTO}:
  Guaranteed Sequential Trajectory Optimization via Sequential Convex
  Programming,} \emph{IEEE International Conference on Robotics and
  Automation}, 2019, pp. 6741--6747.
\newblock \doi{10.1109/ICRA.2019.8794205}.

\bibitem[{A{\c{c}}ikme{\c{s}}e et~al.(2008)A{\c{c}}ikme{\c{s}}e, Blackmore,
  Scharf, and Wolf}]{Acikmese2008}
A{\c{c}}ikme{\c{s}}e, B., Blackmore, L., Scharf, D.~P., and Wolf, A.,
  \enquote{Enhancements on the Convex Programming Based Powered Descent
  Guidance Algorithm for Mars Landing,} \emph{AIAA/AAS Astrodynamics Specialist
  Conference and Exhibit}, 2008.
\newblock \doi{10.2514/6.2008-6426}.

\bibitem[{NAS(2015)}]{NASA2015}
\enquote{NASA Technology Roadmaps TA 9: Entry, Descent, and Landing Systems,}
  2015.

\bibitem[{Dueri et~al.(2017)Dueri, Acikmese, Scharf, and Harris}]{Dueri2017}
Dueri, D., Acikmese, B., Scharf, D.~P., and Harris, M.~W., \enquote{Customized
  Real-Time Interior-Point Methods for Onboard Powered-Descent Guidance,}
  \emph{JGCD Special Issue on Computational Guidance and Control}, Vol.~40,
  No.~2, 2017.
\newblock \doi{10.2514/1.G001480}.

\bibitem[{Elango et~al.(2022)Elango, Kamath, Yu, Acikmese, Mesbahi, and
  Carson}]{Elango2022}
Elango, P., Kamath, A.~G., Yu, Y., Acikmese, B., Mesbahi, M., and Carson,
  J.~M., \enquote{A Customized First-Order Solver for Real-Time Powered-Descent
  Guidance,} \emph{AIAA SCITECH}, 2022.
\newblock \doi{10.2514/6.2022-0951}.

\bibitem[{Kamath et~al.(2022)Kamath, Elango, Yu, Mceowen, Chari, Carson~III,
  and Açıkmeşe}]{kamath2022realtime}
Kamath, A.~G., Elango, P., Yu, Y., Mceowen, S., Chari, G.~M., Carson~III,
  J.~M., and Açıkmeşe, B., \enquote{Real-Time Sequential Conic Optimization
  for Multi-Phase Rocket Landing Guidance,} \emph{arXiv preprint
  arXiv:2212.00375}, 2022.
\newblock \doi{10.48550/arXiv.2212.00375},
  \urlprefix\url{https://doi.org/10.48550/arXiv.2212.00375}.

\bibitem[{Kamath et~al.(2023)Kamath, Elango, Kim, Mceowen, Yu, Carson, Mesbahi,
  and Acikmese}]{kamath2023customized}
Kamath, A.~G., Elango, P., Kim, T., Mceowen, S., Yu, Y., Carson, J.~M.,
  Mesbahi, M., and Acikmese, B., \enquote{Customized Real-Time First-Order
  Methods for Onboard Dual Quaternion-based 6-DoF Powered-Descent Guidance,}
  \emph{AIAA 2023-2003, Session: Entry, Descent and Landing GN\&C Technology
  V}, AIAA, 2023.
\newblock \doi{10.2514/6.2023-2003},
  \urlprefix\url{https://doi.org/10.2514/6.2023-2003}.

\bibitem[{Scharf et~al.(2015)Scharf, Ploen, and Acikmese}]{Scharf2015}
Scharf, D.~P., Ploen, S.~R., and Acikmese, B.~A.,
  \enquote{Interpolation-Enhanced Powered Descent Guidance for Onboard Nominal,
  Off-Nominal, and Multi-X Scenarios,} \emph{SciTech}, 2015.
\newblock \doi{10.2514/6.2015-0850}.

\bibitem[{Kim et~al.(2022)Kim, Elango, Malyuta, and Acikmese}]{Kim2022}
Kim, T., Elango, P., Malyuta, D., and Acikmese, B., \enquote{Guided Policy
  Search using Sequential Convex Programming for Initialization of Trajectory
  Optimization Algorithms,} \emph{American Control Conference}, 2022.
\newblock \doi{10.23919/ACC53348.2022.9867151}.

\bibitem[{Li and Gong(2022)}]{Li2022}
Li, W., and Gong, S., \enquote{Free Final-Time Fuel-Optimal Powered Landing
  Guidance Algorithm Combing Lossless Convex Optimization with Deep Neural
  Network Predictor,} \emph{Applied Sciences}, Vol.~12, No.~7, 2022, p. 3383.
\newblock \doi{10.3390/app12073383}.

\bibitem[{Bertsimas and Stellato(2021)}]{Bertsimas2021}
Bertsimas, D., and Stellato, B., \enquote{The Voice of Optimization,}
  \emph{Machine Learning}, Vol. 110, 2021, pp. 249, 277.
\newblock \doi{10.1007/s10994-020-05893-5}.

\bibitem[{Bertsimas and Stellato(2022)}]{Bertsimas2022}
Bertsimas, D., and Stellato, B., \enquote{Online Mized-Integer Optimization in
  Milliseconds,} \emph{Informs Journal on Computing}, Vol.~34, No.~4, 2022, pp.
  2229, 2248.
\newblock \doi{10.1287/ijoc.2022.1181}.

\bibitem[{Brezak et~al.(2012)Brezak, Bacek, Majetic, Kasac, and
  Novakovic}]{Brezak2012}
Brezak, D., Bacek, T., Majetic, D., Kasac, J., and Novakovic, B., \enquote{A
  Comparison of Feed-Forward and Recurrent Neural Networks in Time Series
  Forecasting,} \emph{Computational Intelligence for Financial Engineering and
  Economics}, 2012.
\newblock \doi{10.1109/CIFEr.2012.6327793}.

\bibitem[{Gers et~al.(2001)Gers, Eck, and Schmidhuber}]{Gers2001}
Gers, F., Eck, D., and Schmidhuber, J., \enquote{Applying LSTM to Time Series
  Predictable through Time-Window Approaches,} \emph{Artificial Neural Networks
  - ICANN}, Vol. 2130, 2001.
\newblock \doi{10.1007/3-540-44668-0_93}.

\bibitem[{Vaswani et~al.(2017)Vaswani, Shazeer, Parmar, Uszkoreit, Jones,
  Gomez, Kaiser, and Polosukhin}]{Vaswani2017}
Vaswani, A., Shazeer, N., Parmar, N., Uszkoreit, J., Jones, L., Gomez, A.,
  Kaiser, L., and Polosukhin, I., \enquote{Attention Is All You Need,}
  \emph{Conference on Neural Information Processing Systems}, Vol.~30, 2017.

\bibitem[{Brown et~al.(2020)Brown, Mann, Ryder, Subbiah, Kaplan, Dhariwal,
  Neelakantan, Shyam, Sastry, Askell, Agarwal, Herbert-Voss, Krueger, Henighan,
  Child, Ramesh, Ziegler, Wu, Winter, Hesse, Chen, Sigler, Litwin, Gray, Chess,
  Clark, Berner, McCandlish, Radford, Sustskever, and Amodei}]{Brown2020}
Brown, T., Mann, B., Ryder, N., Subbiah, M., Kaplan, J., Dhariwal, P.,
  Neelakantan, A., Shyam, P., Sastry, G., Askell, A., Agarwal, S.,
  Herbert-Voss, A., Krueger, G., Henighan, T., Child, R., Ramesh, A., Ziegler,
  D., Wu, J., Winter, C., Hesse, C., Chen, M., Sigler, E., Litwin, M., Gray,
  S., Chess, B., Clark, J., Berner, C., McCandlish, S., Radford, A.,
  Sustskever, I., and Amodei, D., \enquote{Language Models are Few-Shot
  Learners,} \emph{Advances in Neural Information Processing Systems}, Vol.~33,
  2020.

\bibitem[{Devlin et~al.(2019)Devlin, Chang, Lee, and Toutanova}]{Devlin2019}
Devlin, J., Chang, M., Lee, K., and Toutanova, K., \enquote{BERT: Pre-training
  of Deep Bidirectional Transformers for Language Understanding,}
  \emph{NAACL-HLT}, Vol.~1, 2019, pp. 4171, 4186.

\bibitem[{Nie et~al.(2023)Nie, Nguyen, Sinthong, and Kalagnanam}]{Nie2023}
Nie, Y., Nguyen, N., Sinthong, P., and Kalagnanam, J., \enquote{A Time Series
  is Worth 64 Words: Long-Term Forecasting with Transformers,} \emph{ICLR},
  2023.

\bibitem[{Briden et~al.(2023)Briden, Siew, Rodriguez-Fernandez, and
  Linares}]{briden2023transformer}
Briden, J., Siew, P.~M., Rodriguez-Fernandez, V., and Linares, R.,
  \enquote{Transformer-based Atmospheric Density Forecasting,} 2023.
\newblock \doi{10.48550/arXiv.2310.16912},
  \urlprefix\url{https://amostech.com/TechnicalPapers/2023/Atmospherics_Space-Weather/Briden.pdf}.

\bibitem[{Calafiore(2010)}]{Calafiore}
Calafiore, G.~C., \enquote{Random Convex Programs,} \emph{Journal on
  Optimization}, Vol.~20, No.~6, 2010.
\newblock \doi{10.1137/090773490}.

\bibitem[{Paszke et~al.(2019)Paszke, Gross, Massa, Lerer, Bradbury, Chanan,
  Killeen, Lin, Gimelshein, Antiga, Desmaison, Kopf, Yang, DeVito, Raison,
  Tejani, Chilamkurthy, Steiner, Fang, Bai, and Chintala}]{NEURIPS2019_9015}
Paszke, A., Gross, S., Massa, F., Lerer, A., Bradbury, J., Chanan, G., Killeen,
  T., Lin, Z., Gimelshein, N., Antiga, L., Desmaison, A., Kopf, A., Yang, E.,
  DeVito, Z., Raison, M., Tejani, A., Chilamkurthy, S., Steiner, B., Fang, L.,
  Bai, J., and Chintala, S., \enquote{PyTorch: An Imperative Style,
  High-Performance Deep Learning Library,} \emph{Advances in Neural Information
  Processing Systems 32}, Curran Associates, Inc., 2019, pp. 8024--8035.
\newblock
  \urlprefix\url{http://papers.neurips.cc/paper/9015-pytorch-an-imperative-style-high-performance-deep-learning-library.pdf}.

\bibitem[{{van} {der}~Maaten and Hinton(2008)}]{Maaten2008}
{van} {der}~Maaten, L., and Hinton, G., \enquote{Visualizing Data using t-SNE,}
  \emph{Journal of Machine Learning Research}, Vol.~9, No. 2605, 2008, pp.
  2579--2605.

\bibitem[{Pedregosa et~al.(2011)Pedregosa, Varoquaux, Gramfort, Michel,
  Thirion, Grisel, Blondel, Prettenhofer, Weiss, Dubourg, Vanderplas, Passos,
  Cournapeau, Brucher, Perrot, and Duchesnay}]{scikit-learn}
Pedregosa, F., Varoquaux, G., Gramfort, A., Michel, V., Thirion, B., Grisel,
  O., Blondel, M., Prettenhofer, P., Weiss, R., Dubourg, V., Vanderplas, J.,
  Passos, A., Cournapeau, D., Brucher, M., Perrot, M., and Duchesnay, E.,
  \enquote{Scikit-learn: Machine Learning in {P}ython,} \emph{Journal of
  Machine Learning Research}, Vol.~12, 2011, pp. 2825--2830.

\end{thebibliography}

\end{document}